\documentclass[10pt,a4paper,fleqn]{amsart}
\usepackage{amssymb,amsmath,latexsym}
\usepackage[varg]{pxfonts}
\usepackage[english]{babel}
\def\constr#1^#2{\mathrel{\mathop{\kern 0pt#1}\limits^{#2}}}
\def\build#1_#2{\mathrel{\mathop{\kern 0pt#1}\limits_{#2}}}

\theoremstyle{plain}

\newtheorem{Theorem}{Theorem}[section]

\theoremstyle{Corollary}
\newtheorem{Corollary}{Corollary}[section]

\newtheorem{Lemma}{Lemma}[section]
\theoremstyle{Proposition}
\newtheorem{Proposition}{Proposition}[section]
\theoremstyle{definition}
\newtheorem{Definition}{Definition}[section]

\newtheorem{remark}{Remark}[section]

\theoremstyle{Abstract}

\numberwithin{equation}{section}

\newcommand{\hf}{\hfill $\diamondsuit$ }    % ajouter le petit losange \`{a} la fin
\newenvironment{pf}{\noindent{\bf Proof.}\enspace}{%\rule{2mm}{2mm}
\hfill Q.E.D.}%\medskip}

\begin{document}
%\title{ New fixed point theorem in generalized Banach Algebra under weak topologie features and application }
%\centerline{\bf \Large  New fixed point theorem in generalized Banach Algebras}
%\centerline{\bf \Large  under weak topologie features and application}
\title[]{\large Generalized form of fixed-point theorems in generalized Banach algebra relative to the weak topology with an application}
%\centerline{\bf Some hybrid fixed point theorems in generalized Banach algebra under weak topologie  with an application}
%\centerline{ \bf   }
\maketitle

\par\vskip0.1cm

%%%%%%%%%%%%%%%%%%%%%%%%%%%%%%%%%%%%%%%%%%%%%%%%%%%%%%%%%%%%%%%%%%%%%%%%%%%%%%%%%%%%%%%%%%%%%%%%%

\author\centerline{{ Aref Jeribi $^{(1)(i)}$,   Najib Kaddachi $^{(1)(ii)}$  and  Zahra Laouar$^{(1)(i)}$}}

\vspace{3mm}

\begin{center}
\emph{$^{(1)}$ Department of Mathematics. University of Sfax.
Faculty of Sciences of Sfax. \\ Soukra Road Km 3.5 B.P. 1171, 3000, Sfax,
Tunisia}.\\

\bigskip

e-mail : \  $^{(i)}$ Aref.Jeribi$@$fss.rnu.tn\ \ \ \ \  $^{(ii)}$ najibkadachi$@$gmail.com

\end{center}
%\centerline{\bf \Large under weak topologie}
\vskip0.2cm
{\footnotesize
\noindent
\parbox{5.5 in}{{\bf\small Abstract.} In this paper, a general hybrid fixed point theorem for the contractive  mappings in generalized Banach spaces is proved via  measure of weak  non-compactness and it is further applied to  fractional integral equations   for proving the existence results for the  solutions under mixed Lipschitz and weakly sequentially continuous conditions. Finally, an example is given to illustrate the result.
  %The main tool used in our considerations is the fixed point theorems for the sum and the product of nonlinear weakly sequentially continuous operators defined on a nonempty, closed and convex subset of generalized Banach algebra.
% we prove fixed point theorems for the sum and the product of nonlinear weakly sequentially continuous operators acting on a generalized Banach algebras.
 % The obtained results are then applied to an fractional integral equations.
  {\small \sloppy{}}}

\noindent
\parbox{5.5 in}{{\bf\small Keywords:} Banach space,  weakly compact, Measure of weak
non-compactness, Fixed point theory, Integral equations.
 {\small \sloppy{ }}}

\bigskip

}

\section{Introduction}
Many nonlinear problems involve the study of nonlinear equations of the form
%In this paper, we are concerned with existence of solutions of the following equation:
\begin{eqnarray} \label{abc}
x=Ax \cdot Bx+Cx, \ \ x\in S
\end{eqnarray}
where $S$ is a nonempty, closed, and convex subset of a Banach algebra $X,$ see for example \cite{1,2,3,4,5} and the references therein. These studies were mainly based on  the Schauder fixed point theorem, and properties of the operators $A, B$ and $C$ (cf. completely
continuous, k-set contractive, condensing, and the potential tool of the axiomatic measure of noncompactness,$\ldots).$ Various attempts have been made in the literature to extend this study to a weak topology, see for example \cite{sofiane,so,wc}. In recent years, the authors Ben Amar et al. have established in \cite{sofiane}  some fixed point theorems for the operator equation $(\ref{abc})$ in Banach algebra under the weak topology, their results were based on the class of weak sequential continuity, weakly compact and weakly condensing conditions. In other direction, the classical Banach contraction principle was extended for contractive maps on spaces endowed with vector-valued metrics by Perov  \cite{{perov1},{perov2}}.  Schauder's fixed point theorem has been extended from Banach spaces   to generalized Banach spaces by Viorel \cite{Viorel} and Krasnosel'skii's  fixed point theorem  has been extended by  Petru and ouahab in \cite{{Ouahab},{Petru}}. More recently, the authors Nieton et al.  \cite{algebra} have
also established some new variants of fixed-point theorems for operator equation $(\ref{abc})$ in vector-valued metrics endowed with an internal composition
law $(\cdot).$ These studies were mainly based on the convexity and the
closure of the bounded domain and properties of the operators $A, B$ and $C.$ Their analysis was carried out via arguments of strong topology.\\
%relative to the strong topology.\\

% Our results extend and improve well-known results in \cite{sofiane,so,algebra}
\par This paper is centered around the following question: under which conditions on its entries, the operator equation $(\ref{abc})$ acting on a generalized Banach algebras with respect to the weak topology, has a solution?\\
% The main goal of this research is to extend some new variants of fixed point theorems for the   in   generalized Banach algebras under weak topology in the spirit of the papers .
Our main results are applied to investigate the existence of solutions for the following coupled system  of quadratic integral
equations of fractional order
\begin{equation}\label{Sy3}
x_i(t)= f_i(t,x_1(t),x_2(t))\cdot \displaystyle\int_{0}^{t}\displaystyle\frac{(t-s)^{\alpha_i-1}}{\Gamma(\alpha_i)}g_i(s,x_1(s),x_2(s))ds+\displaystyle\sum _{k=1}^{m}I^{\beta_{i}^{k}}h_{i}^{k}(t,x_1(t),x_2(t)),\ \ i=1,2$$
\end{equation}
%$$\left\{
%  \begin{array}{ll}
%  x(t)= f_1(t,x(t),y(t))\cdot \displaystyle\int_{0}^{t}\displaystyle\frac{(t-s)^{\alpha_1-1}}{\Gamma(\alpha_1)}g_1(s,x(s),y(s))ds+\displaystyle\sum _{k=1}^{m}I^{\beta_{1}^{k}}h_{1}^{k}(t,x(t),y(t)),\ \ \hfill{ t \in J,} & \hbox{} \\ \\
% y(t)= f_2(t,x(t),y(t)) \cdot \displaystyle\int_{0}^{t}\displaystyle\frac{(t-s)^{\alpha_2-1}}{\Gamma(\alpha_2)}g_2(s,x(s),y(s))ds+\displaystyle\sum _{k=1}^{m}I^{\beta_{2}^{k}}h_{2}^{k}(t,x(t),y(t)),\ \ \hfill{ t \in J,} & \hbox{}
%\end{array}
%\right.
%$$
 where   $k \in \{1,...,m\},$  $ t,\alpha_{i}\in (0,1),$    $I^{\beta_{i}^{k}}$ is the fractional Pettis integral  of order $\beta_{i}^{k}>0,$  $\Gamma(\cdot)$ is the Gamma function and the functions  $f_i,g_i,h_i^k$ are given functions,
whereas $x_1=x_1(t)$ and $x_2=x_2(t)$  are unknown functions.\\

\par The present paper is arranged as follows. The next section, we give some preliminary results for future use. Moreover,  we shall extend the result of B. Amar, Jeribi and  Mnif [\cite{compact}, Theorem $2.5$]. In addition, we establish the fixed point theorem by using the concept of the measure of weak non compactness in generalized Banach space $($see Theorem \ref{condensing}$)$. Note our result $($Theorem \ref{condensing}$)$ improves and generalizes Theorem $3.2$ in \cite{ther3.2}. Apart from that we introduce a class of  generalized Banach algebras satisfying certain sequential conditions called here the condition $(\mathcal{GP})$ $($see Definition \ref{defi}   $ )$. In section $3,$ we present  a  collection of new fixed point theorems in generalized Banach algebras satisfying condition $(\mathcal{GP})$. Our results extend and improve well-known results in \cite{sofiane,so,algebra}. In the last section of this manuscript, we apply Theorem $\ref{3}$ to discuss the existence of solutions for a  system of fractional integral equations $(\ref{Sy3})$ and  an example is given to explain the applicability of the results.

\section{ Preliminaries and Results }
\noindent Let $(X,\|\cdot\|)$ be a generalized Banach space in the sense of Perov such that the vector-valued norm $\|\cdot\|:X \longrightarrow \mathbb{R}^n_+$ is given by
 $$\|x\|=\left(%
\begin{array}{ccc}
  \|x\|_1  \\
  \vdots \\
   \|x\|_n \\
\end{array}%
\right),\ \  x \in X $$
with $\|\cdot\|_i, i=1,...,n$ define $n$ norm on $X.$ We denote by  $\theta$ the zero element of $X$ and ${B}_r =B(\theta,r) $ the  closed  ball  centered at $\theta$ with radius $r \in \mathbb{R}^{n}_{+}.$ For more details, the reader may consult the monograph of John R. Graef, Johnny Henderson and Abdelghani Ouahab \cite{monograph}.\\
Let  $ (Y_i,\tau_i)_{i \in I} $ be a family of topological spaces and   let $f_i : X \longrightarrow Y_i,\ \ i \in  I$ be a linear and continuous  mappings. We define the smallest  topology on $X$ such that all the mappings $f_i$ remain continuous with respect to this topology. Its basic open sets are of the form $ \displaystyle\cap_{i \in J}f_i^{-1}u_i, $ with $J$ a finite subset of $I$ and $u_i \in  \tau_i$ for each $i.$ This topology is called the weak topology on $X$ generated by the $(f_i)_{i \in I} $ and we denote it by $\sigma(X,(f_i)_{i \in I} ).$
So that a sequence $(x_n)_n$ in $X$ converges to $x$ in $\sigma(X, (f_i)_{i \in I} )$  if and only if $(f_i(x_n))_n$ converges to $f_i ( x ),$ for all $i \in I.$ In generalized Banach spaces,   weakly open subset, weakly closed subset and weak compactness, are similar to those for usual  Banach spaces. We denote by $\mathcal{B}(X)$  the collection of all nonempty bounded subsets of X and $\mathcal{W}(X)$  is a subset of $\mathcal{B}(X)$ consisting of all weakly compact subset of $X.$ If $x,y \in \mathbb{R}^{n},$ $x=(x_1,...,x_n),\ \ y=(y_1,...,y_n)$ by $x\leq y$ we means $x_i \leq y_i$ for all $i=1,...,n.$
%If $c \in \mathbb{R},$ then $x\leq c$ means $x_i \leq c$ for each $i=1,...,n.$

\bigskip
\begin{Definition}
A square matrix of real numbers $M$ is said to be convergent to zero if and only if its spectral radius $\rho(M) $ is strictly less than $1.$ In other words, $|\lambda|<1,$ for every $\lambda \in \mathbb{C}, $ with $det(M-\lambda I)=0,$ where $I$ denotes the unit matrix of $\mathcal{M}_{n \times n}(\mathbb{R}_+).$ \hf
\end{Definition}
\bigskip

\begin{Lemma}\label{lem}\cite{matric} Let $M $ The following assertions are equivalent:\\
\noindent{$(i)$}  $M$ is a matrix convergent to zero,\\
\noindent{$(ii)$}  $M^k \rightarrow 0 $ as $k \rightarrow \infty,$\\
\noindent {$ (iii)$}  The matrix $(I-M)$ is nonsingular and $$(I-M)^{-1}= I+M+M^{2}+...+M^k+...,$$\\
\noindent {$(iv)$} The matrix $(I-M)$ is nonsingular and $(I-M)^{-1}$ has nonnegative elements.\hf
\end{Lemma}

\bigskip

\begin{Definition}\cite{rus}
Let $X$ be a generalized Banach space. An operator $ T:X\longrightarrow X$ is said to be contractive if there exists a matrix $M \in  \mathcal{M}_{n \times n}(\mathbb{R}_+)$ convergent to zero such that
$$\|Tx-Ty\|\leq M\|x-y\|,\ \ \text{ for all }  x, y \\  \text { in } X. $$
For $n = 1$, we recover the classical Banach's contraction fixed point result. \hf
\end{Definition}

\bigskip

\begin{Definition}
 An operator $T:X\longrightarrow X$ is called weakly sequentially continuous on $X,$ if for every sequence $(x_n)_{n \in \mathbb{N}}$ with
 %\vskip0.2cm\hskip4cm
   $x_n \rightharpoonup x \ \ \text{  we have  }\ \  Tx_n \rightharpoonup Tx,  $   here $\rightharpoonup$  denotes weak convergence. \hf
   \end{Definition}

\bigskip

\begin{Definition}
 An operator $T:X\longrightarrow X$ is said to be weakly compact if $T(V)$ is relatively weakly compact for every bounded subset $ V\subseteq X. $ \hf
\end{Definition}

\bigskip
\begin{Theorem}\cite{Viorel,sch2}\label{C1}
Let $X$ be a generalized Banach space, $S$ be a nonempty, compact and convex subset of $X$ and $A: S\longrightarrow S$ be a continuous operator. Then $A$ has at least a fixed point in $S.$ \hf
\end{Theorem}
\bigskip
%\begin{Theorem} \label{the1}
%Let $X$ be a generalized Banach space, $S$ be a nonempty, weakly compact, convex and  $A:S\longrightarrow S$ be a sequentially weakly continuous operator. Then $A$ has, at least, a fixed %point in $S.$ \hf
%\end{Theorem}
%\begin{pf}
%It suffices to prove that $T$ is weakly continuous, so that the Schauder fixed point theorem applies. Now for each weakly closed subset $M$ of $S$,
%$T^{-1}(M)$ is sequentially closed in $S$, hence weakly compact by the Eberlein-\u{S}mulian's
%theorem \cite{} and $T^{-1}(M)$ is weakly closed. Hence $T$ is weakly continuous.
%\end{pf}
%
\bigskip
\noindent Now, we  will prove the following  theorem.
\begin{Theorem} \label{1}
Let $S$ be a nonempty, closed and convex subset of a generalized Banach space $X.$ Let $T : S \longrightarrow S$ be a weakly  sequentially continuous mapping. If $T(S)$ is relatively weakly
compact, then $T$ has a fixed point. \hf
\end{Theorem}
\begin{pf}Let $M=\overline{co}(T(S))$ the closed convex hull of $T(S).$ Because $T(S)$ is relatively weakly compact, then $M$ is weakly compact convex subset of $X.$ On the other hand, $T(M)\subset T(S) \subset \overline{co}(T(S)=M $ ie : $T$ maps $M$  into itself. Since $T$ is weakly sequentially continuous on $M,$  we infer that $T$ is weakly continue on $M.$ (consider $X=(X,\sigma(X,(f_i)_{i \in {I}}))$ the space endowed with the weak topology and note that $T:M \longrightarrow M $ is continuous with M is compact). The use of Theorem \ref{C1}
 achieves the proof
\end{pf}

\bigskip
\begin{Definition}
the  measure of weak non-compactness of the generalized Banach space $X$ is the map $\mu: \mathcal{B}{(X)}\longrightarrow \mathbb{R}_{+}^{n}$ defined in the following way
$$ \mu(S)= \inf \big\{ r \in \mathbb{R}^{n}_{+},\ \ \text{ there exists } K\in \mathcal{W}(X)\text{ such that } S \subseteq K+ {B}_r \big  \},$$
 for all $S \in \mathcal{B}(X).$ \hf
\end{Definition}

 \begin{Lemma} Let $S_1$, $S_2$ be two elements of $\mathcal{B}(X).$ Then the functional $\mu$  has the properties:\\
 \noindent $(i)$ $S_1 \subseteq S_2$ implies $\mu(S_1)\leq \mu(S_2),$\\
 \noindent $(ii)$ $\mu(\overline{S_1^w})=\mu(S_1),$\\
 \noindent $(iii)$ $\mu(S_1)=0,$ if and only if, $ \overline{S_1^w} \in \mathcal{W}{(X)},$ where $\overline{S_1^w}$ is the weak closure of $S_1,$\\
\noindent $(iv)$  $\mu(co(S_1))=\mu(S_1),$ where $co(S_1)$ is the convex hull of $S_1,$ \\
\noindent $(v)$  $\mu(S_1+S_2)\leq \mu(S_1)+\mu(S_2), \mu(\{a\}+S_1)\leq \mu(S_1),$\\
\noindent $(vi)$ $\mu(S_1 \cup S_2)= \max \{\mu(S_1),\mu(S_2)\},$\\
\noindent $(vii)$ $\mu(\lambda S_1)= |\lambda|\mu(S_1),$ for all $\lambda \in \mathbb{R}.$\\
 \hf
 \end{Lemma}
\begin{pf}
The statements $(i)$ and $(iv)$ - $(viii)$ are immediate consequences of the definition of $\mu$
.  Let us prove $(ii)$. From the definition of $\mu$  there exists a subset $K_1 \in \mathcal{W}(X)$ and $r_1 \in \mathbb{R}^{n}_+\setminus\{0_{\mathbb{R}^{n}}\}$  such that    $S_1  \subseteq  K_1+ {B}(\theta,r_1+\mu(S_1)),$ then $   S_1  \subseteq  \overline{co}{K_1}+ {B}(\theta,r_1+\mu(S)).$   By the Kerein-\u{S}mulian theorem   $\overline{co}{K_1}$ is weakly compact. Since  $ {B}(\theta,r_1+\mu(S_1))$ is weakly closed, we infer that
  $    \overline{co}{K_1}+ {B}(\theta,r_1+\mu(S_1))$  is weakly closed. So  $ \overline{S}_1^w \subseteq \overline{co}{K_1}+ {B}(\theta,r_1+\mu(S_1))$ implies that
  $\mu(\overline{S}_1^w )\leq r_1+\mu(S_1) .$ Letting $r_1 \to 0_{\mathbb{R}^{n}}$ in the above inequality, we get $\mu(\overline{S}_1^w )\leq \mu(S_1) .$ The reverse inequality is obvious  .
  The proof of the "only if" part of $(iii).$   By the definition of $\mu,$  there exists a subset $K_1 \in \mathcal{W}(X)$ and $r_1 \in \mathbb{R}^{n}_+\setminus\{0_{\mathbb{R}^{n}}\}$  such that
$$\overline{S}_1^w \subseteq \overline{K}_1^w+{B}(\theta, r_1). $$
Let $\{x_n\}_{n\in \mathbb{N}}$ be a sequence on $\overline{S}_1^w.$  So there existe  sebsequences $\{y_n\}_{n \in \mathbb{N}}$ and $\{z_n\}_{n \in \mathbb{N}}$ of $\overline{K}_1^w$ and  ${B}(\theta, r_1)$ respectively, such that
$x_n=y_n+z_n.$
 If $r_1  \longrightarrow 0_{\mathbb{R}^{n}},$ we deduce that  $\{x_n\}$ has a weakly convergence subsequence, So $\overline{S}_1^w$ is weakly compact, while the "if" part is trivial.\\

\ \ \ \

\end{pf}
\bigskip

\noindent For $n = 1$, we recover the classical De Blasi measure of weak non-compactness.

\bigskip

\begin{Definition} Let $S$ be a nonempty subset of a generalized Banach space $X$ and $M \in \mathcal{M}_{n \times n}(\mathbb{R}_+) $ is a matrix convergent to zero. If $T$ maps $S$ into $X,$ we say that \\
\noindent $(i)$ $T$ is $M$-set contraction with respect to $ \mu$ if $T$ is bounded and for any bounded subset $V  \text{  in  } S$ and  $\mu(T(V))\leq M\mu(V), $\\
\noindent $(ii)$ $T$ is condensing with $\mu$ if $T$ is bounded and $\mu(T(V))<\mu(V)$ for all bounded subsets $V$ of $S$ with $\mu(V)>0_{\mathbb{R}^{n}}.$\hf
\begin{remark}
If we assume that $T$ is $M$-set contraction, then $T$ is condensing. Indeed, let $V$ be a bounded subset on $S$ with $\mu(V)>0_{\mathbb{R}^{n}}.$ We claim that $M\mu(V)<\mu(V).$ If not we obtain $(I-M)\mu(V)\leq 0_{\mathbb{R}^{n}}.$ Since   $(I-M)^{-1}$ has nonnegative elements , it follows that  $\mu(V)=0_{\mathbb{R}^{n}}.$ which is a contradiction with $\mu(V)>0_{\mathbb{R}^{n}}.$
\hf
\end{remark}

\bigskip
 \noindent By using  the concept of a  measure  of weak noncompactness in vector-valued Banach spaces, we obtain  the following fixed point theorem.
\begin{Theorem} \label{condensing}   Let  $S$ be a nonempty closed convex subset of a generalized Banach space  $X.$ Let $T:S\longrightarrow S$ be a weakly seqeuntially  continuous operator and condensing with respect to $\mu .$   If $T(S)$ is bounded, then $T$ has a fixed point in $S.$ \hf
\end{Theorem}
\begin{pf}
let $x_0$ be fixed in $S$ and define the set
$$\Sigma= \{K:\ \ x_0 \in K \subseteq S,\ \ K \text{   is closed, convex, bonded and   } T(K) \subseteq K \}. $$
Clearly, $\Sigma \neq \emptyset$ since  $\overline{co}(T(S)\cup \{x_0\})\subseteq S$ and  we have
$$T(\overline{co}(T(S)\cup \{x_0\})\subseteq T(S) \subseteq \overline{co}(T(S)\cup \{x_0\}). $$
 Which shows that $\overline{co}(T(S)\cup \{x_0\})\in \Sigma.$
 If we consider $M=\displaystyle\cap_{K \in \Sigma }K,$ then $x_0 \in M \subseteq S, M$ is also a closed convex subset and $T(M)  \subseteq M.$ Therefore, we have that $ M \in \Sigma.$  We will prove that $M$ is weakly compact. Denoting by $M_0=\overline{co}(T(M)\cup{x_0})$, we have $ M_0 \subseteq M $, which implies that $T(M_0)\subseteq T(M) \subseteq M_0$. Therefore $M_0 \in \Sigma $, $M \subseteq M_0 $. Hence $M=M_0.$ Since $M$ is weakly closed, it suffices to show that $M$ is relatively weakly compact. If $\mu(M)>0_{\mathbb{R}^{n}},$ we obtain
$$\mu(M)=\mu(\overline{co}(T(M)\cup \{x_0\}))< \mu(M) $$
which is a contradiction. Hence, $\mu(M)=0_{\mathbb{R}^{n}}$ and so $M$ is relatively weakly compact. Now, $T$ is weakly continuous on $M.$ Consider $X=(X,\sigma(X,(f_i)_{i \in {I}}))$   the space endowed with the weak topology.  Hence, an application of Theorem \ref{C1}  shows that $T$ has at least one fixed point in $M.$
\end{pf}

% An operator $T:S\longrightarrow X$ is said to be $M$-$\mu$-contraction if $T$ is bounded and for any bounded subset $V  \text{  in  } S,$  $\mu(T(V))\leq M\mu(V). $\hf
\end{Definition}

\bigskip

\begin{Definition}\cite{algebra} A generalized normed algebra $X$ is an algebra endowed with a norm
satisfying the following property
$$\text{  for all  }  x, y \in X \ \ \|x.y\|\leq \|x\|\|y\|,$$ where $$\|x.y\|=\left(%
\begin{array}{ccc}
  \|x.y\|_1  \\
  \vdots \\
   \|x.y\|_n \\
\end{array}%
\right)$$ and $$\|x\|\|y\|=\left(%
\begin{array}{ccc}
  \|x\|_1\|y\|_1  \\
  \vdots \\
   \|x\|_n\|y\|_n \\
\end{array}%
\right).$$ A complete generalized normed algebra is called a generalized Banach algebra.\hf

\end{Definition}

\bigskip
\noindent Because the product of two sequentially weakly continuous functions in generalized Banach algaebra is not necessarily sequentially weakly continuous, we will introduce:
\begin{Definition}\label{defi}
we will say that the generalized Banach algebra $X$ satisfies conditions $(\mathcal{G P})$\ \ if
 \vskip0.2cm\hskip0.2cm $(\mathcal{G P})\left\{
 \begin{array}{ll}
 \mbox{~For any sequences ~}\{x_n\} \mbox{~and~} \{y_n\} \mbox{~of~} X \mbox{~such that~} x_n\rightharpoonup x \mbox{~and~} y_n \rightharpoonup y,&\hbox{}\\
\mbox{~then~} x_n \cdot y_n \rightharpoonup x\cdot y; \mbox{~here~} \rightharpoonup \mbox{~denotes weak convergence~}.&\hbox{}
\end{array}
 \right.$ \hf
\end{Definition}
\noindent If $X$ is Banach algebra, we recover the classical sequential condition $(\mathcal{P})$ [\cite{so},Definition $3.1$]
\bigskip
%\begin{Definition}
%Let $f:[0,T]\rightarrow X $ be a function. The fractional Pettis integral of the
%function $f$ of order $ \alpha \in  \mathbb{R_+}$  is defined by
%$$ I^{\alpha}f(t)=\int_{0}^{t}\frac{(t-s)^{\alpha}}{\Gamma(\alpha)}f(s)ds,  $$
%where the sign $ "\int " $  denotes the Pettis integral and $\Gamma$ is the gamma function.\hf
%\end{Definition}

\begin{Lemma}\label{lemma}
If $K$ and $ K'$ are weakly compact subsets of generalized Banach algebra $X$ satisfying the condition $(\mathcal{G P} ),$ then $K \cdot K'$ is weakly compact.\hf
\end{Lemma}

\bigskip
\begin{pf}
We will show that $K\cdot K'$ is weakly sequentially  compact. For that, let $\{z_n\}_n$  be any sequence of $K \cdot K '.$ So, there exist  sequences $\{x_n\}_n$ and $\{x'_n\}_n$ of $K$ and $K'$
respectively.
% $\{x_n\}_n$ be any sequence of $K$ and let $\{x'_n\}_n$ be any sequence of $K'$.
 By hypothesis, there is a renamed subsequence $\{x_n\}_n$ such that $ x_n \rightharpoonup x \in K.$ Again, there is a renamed subsequence $\{x'_n\}_n$  of $K'$ such that $ x'_n \rightharpoonup x' \in K.$  This, together with condition $(\mathcal{GP})$ yields that
$$z_n\rightharpoonup z= x\cdot x'.$$
This implies that $K  \cdot K'$ is weakly sequentially  compact. Hence, an application of the Eberlein-\u{S}mulian's theorem yields that $K \cdot K' $ is weakly compact.
%The proof of Lemma is checked in a similar way to that in Lemma $3.1$ in \cite{sofiane}
\end{pf}
\bigskip
%\noindent Reasoning as in the proof of theorem $2.5$ in \cite{condensing}, we get the following result
\bigskip

\section{ Existence solutions of Equation (\ref{abc})}

\noindent In this section,  we are prepared to state our first fixed point theorems in  generalized Banach algebra in order to found the existence of solutions for the operator equation $( \ref{abc}  )$ under weak topology  in the case where $A,B$ and $C$ are weakly sequentially continuous operators.
%First of our main results is the following.
\begin{Theorem} \label{3.1} Let $S$ be a nonempty, closed, and convex subset of a generalized Banach algebra $X$ satisfying the sequential condition  $(\mathcal{GP}).$ Assume that  $A, C:X \longrightarrow X$ and $B:S \longrightarrow X$ are three weakly sequentially  continuous operators satisfying the following conditions:\\
\noindent $(i)$ $A$ and $C$ are contractive with  Lipschitz matrix $M_A$ and $M_C$ respectively,\\
\noindent $(ii)$ $A$ is regular on $X,$ $($i.e. $A$ maps $X$ into the set of all invertible elements of $X),$\\
\noindent $(iii)$ $B(S)$ is bounded,\\
\noindent $(iv)$ $x=Ax\cdot By+Cx,\ \ y \in S \Rightarrow  x \in S,$ and \\
\noindent $(v)$ $\left(\displaystyle\dfrac{I-C}{A}\right)^{-1}$ is weakly compact.\\
Then, the operator equation $($\ref{abc}$)$ has, at least, a solution   in $S$ as soon as $\|B(S)\|M_A+M_C \in \mathcal{M}_{n \times n}(\mathbb{R}_+) $ is a matrix converging to zero.\hf
\bigskip
\begin{remark}
If $M_A= (a_{ij})_{1\leq i,j\leq n}, \ \ M_C= (\overline{a}_{ij})_{1\leq i,j\leq n}  $ and $\|B(S)\|=(b_j)_{1\leq j \leq n}$ then
 $$\|B(S)\|M_A+M_C=\left(
\begin{array}{ccc}
  b_1a_{11}+\overline{a}_{11}  & \ldots & b_na_{1n}+ \overline{a}_{1n} \\
  \vdots & \ddots & \vdots  \\
  b_1a_{n1}+\overline{a}_{n1}  &\ldots & b_na_{nn}+\overline{a}_{nn}, \\
\end{array}
\right).$$ \hf
\end{remark}
\bigskip
 %where $bM_A=(b_ja_{ij})_{1\leq i,j\leq n}$ and $M_C=(\overline{a}_{ij})_{1\leq i,j\leq n}.$
\end{Theorem}

\noindent \textbf{ Proof of Theorem \ref{3.1}. } It is easy to check that the vector $x\in S$ is a solution for the operator equation $x=Ax\cdot Bx+Cx,$  if and only
if $x$ is a fixed point for the inverse operator   $T :=\left(\displaystyle\dfrac{I-C}{A}\right)^{-1}B.$ The use of assumption $(i)$ as well as perov's theorem  \cite{banach gen} leads to for each $y\in S,$ there is a unique $x_y\in X$ such that $x_y=Ty.$ Indeed,
let $y\in S$ be arbitrary and let us define the mapping
$\varphi_y: X  \longrightarrow X$ by the formula $$\varphi_y(x)=Ax \cdot By+Cx.$$
Notice that $\varphi_y$ is  contractive  with a Lipschitz matrix $\|B(S)\|M_A+M_C.$ Applying Perov's theorem  \cite{banach gen}, we obtain that $\varphi_y$ has a unique fixed point  in $S,$ say $x_y.$ From assumption $(ii),$ it follows that  the  operator $\left(\displaystyle\dfrac{I-C}{A}\right)^{-1}$ is well defined on $B(S).$\\
%Let $x_1,x_2 \in X,$ the use of the assumption $(i)$ reveals that
%\begin{eqnarray*}
%\|\varphi_y(x_1)-\varphi_y(x_2)\| &\leq & \|Ax_1By-Ax_2By\|+\|Cx_1-Cx_2\| \\
 %                            &\leq & \|Ax_1-Ax_2\|\|By\|+ \|Cx_1-Cx_2\| \\
 %                           &\leq & M_A\|x_1-x_2\|\|By\|+M_C\|x_1-x_2\|\\
  %                           &\leq & (bM_A+M_C) \|x_1-x_2\|.
%\end{eqnarray*}
%Taking into occount the assumption $(iv),$ together with the ,we get the operator $T :=\Bigg(\displaystyle\dfrac{I-C}{A}\Bigg)^{-1}B: S\longrightarrow S$ is well defined.
%%By using Theorem $10.1$ in \cite{banach gen} there is a unique point $x\in X$ such that $x=\varphi_y(x),$
%% hence, the operator $T: S\longrightarrow X$ is well defined. Notice also by assumption $(iii)$ we have $T(S)\subset S.$
% Moreover, $T(S)$ is relatively weakly compact. In order to achieve
%the proof, we will apply Theorem \ref{1} , Hence, we only have to prove that the operator $T$ is
%    % Now, we show that $\Bigg(\displaystyle\dfrac{I-C}{A}\Bigg)^{-1}B$ is
%weakly sequentially continuous. To see this,
Let $\{x_n, \ n\in \mathbb{N} \}$ be a weakly convergent sequence of $S$ to a point $x\in S$ and let $y_n=Tx_n.$ Then the relation $y_n=Ay_n\cdot Bx_n+Cy_n$ holds and, therefore $\{y_n, \ n\in \mathbb{N} \}\subset S$ in view of assumption $(iv).$
 Since $T(S)$ is relatively weakly compact, there is a   subsequence $(x_{n_{k}})$ of $\{x_n, \ n\in \mathbb{N} \}$ such that $y_{n_{k}}=Tx_{n_{k}}\rightharpoonup y,$ { for some } $y\in X.$ The weak sequential closedness of $S$ gives $y\in S.$ Making use of the condition $(\mathcal{GP}),$ together with the assumptions on $A,B$ and $C,$ enables
us to have $$y_{n_{k}}=Ay_{n_{k}}\cdot Bx_{n_{k}}+Cy_{n_{k}}\rightharpoonup Ay\cdot Bx+Cy.$$ This implies that $y=Ay\cdot Bx+Cy.$ Consequently $Tx_{n_{k}}\rightharpoonup Tx$ in light of
assumption $(ii).$

\noindent Now a standard argument shows that $Tx_n\rightharpoonup Tx.$
    Suppose the contrary, then there exists  a weakly neighborhood $V^{w}$ of $Tx$ and a subsequence $(x_{n_j})_j$  of  $\{x_{n}, \ n\in \mathbb{N} \}$  such that $Tx_{n_{j}} \notin V^{w},$ for all $j\geq 1.$  Arguing as
before, we may extract a subsequence  $(x_{n_{j_{k}}})_{k \in \mathbb{N}}$ of $\{x_{n_j}, \ j\in \mathbb{N} \}$ verifying  $Tx_{n_{j_{k}}}\rightharpoonup y,$
%However $$ \text{  for all  }  j \in \mathbb{N}, y_{n_{j}} \in T(S).$$
% \Bigg(\displaystyle\dfrac{I-C}{A}\Bigg)^{-1}B(S).  $$
  which is a contradiction and consequently $T$ is weakly sequentially continuous. Hence, $T$ has, at least, one fixed point $x$ in $S$ in view of Theorem \ref{1}.
 \hfill{Q.E.D.}

\bigskip
\noindent An interesting fixed point result of Theorem \ref{3.1} is
\bigskip
\begin{Corollary} Let $S$ be a nonempty, closed, and convex subset of a generalized Banach algebra $X$ satisfying the sequential condition  $(\mathcal{GP}).$ Assume that  $A, C:X \longrightarrow X$ and $B:S \longrightarrow X$ are three weakly sequentially  continuous operators satisfying the following conditions:\\
\noindent $(i)$ $A$ and $C$ are contractive with Lipschitz matrix $M_A$ and $M_C$ respectively,\\
\noindent $(ii)$ $A$ is regular on $X,$\\
%\noindent $(iii)$ $A,B$ and $C$ are \\
\noindent $(iii)$ $A(S),B(S)$ and $C(S)$ are relatively weakly compact,\\
\noindent $(iv)$ $x=Ax{\cdot}By+Cx,\ \ y \in S \Rightarrow  x \in S. $\\
Then, the operator equation $($\ref{abc}$)$ has, at least, a solution   in $S$ as soon as $\|B(S)\|M_A+M_C $ is a matrix converging to zero.\hf

\end{Corollary}
%\vskip0.2cm\hskip4cm$b=\|B(S)\|.$  \hf
%\par Now we brove the following results, which we need in the sequel.

%\noindent Now our attention is directed towards the case where only $B$ is weakly compact in $S.$
\bigskip
\begin{Theorem}\label{3} Let $S$ be a nonempty, bounded, closed, and convex subset of a generalized Banach algebra $X$ satisfying the sequential condition  $(\mathcal{GP}).$ Assume that  $A, C:X \longrightarrow X$ and $B:S \longrightarrow X$ are three weakly sequentially  continuous operators satisfying the following conditions:\\
\noindent $(i)$ $A$ and $C$ are contractive with Lipschitz matrix $M_A$ and $M_C$ respectively,\\
\noindent $(ii)$ $A$ is regular on $X$,\\
\noindent $(iii)$ $B$ is weakly compact, \\
\noindent $(iv)$ $x=Ax\cdot By+Cx,\ \ y \in S \Rightarrow  x \in S. $\\
 Then, the operator equation $($\ref{abc}$)$ has, at least, a solution   in $S$ as soon as $\|B(S)\|M_A+M_C $ is a matrix converging to zero.\hf

\end{Theorem}

\bigskip

\noindent Before proving the theorem, we need to establish two lemmas.
\bigskip
\begin{Lemma} \label{lemma3.2} Let $X$ be  a generalized Banach algebra. If $S\in \mathcal{B}(X)$ and $K\in \mathcal{W}(X),$ then
%For any bounded subset $S$ of  and for any weakly compact subset $K$ of $X,$ we have
%\vskip0.2cm\hskip4cm
$\mu \left(S \cdot K\right)\leq \|K\|\mu (S).$ \hf
\end{Lemma}
\begin{pf} Assume that $\|K\|>0_{\mathbb{R}^{n}}.$
%We note that $K$ is bounded, then $\|K\|$ exists. Next,  $>0,$
% where $\|K\|=\left(%
%\begin{array}{ccc}
% \|K\|_{1}  \\
%  \vdots \\
%  \|K\|_{n} \\
%\end{array}%
%\right).$
   By the definition of $\mu,$  there exists a subset $K'$ of $\mathcal{W}(X)$  and    $ r  \in {\mathbb{R}^{n}_+}\setminus \{0_{\mathbb{R}^{n}}\}$ such that
%$$S\subseteq K'+\overline{\mathcal{B}_i}\left(0,\mu_i(S)+\dfrac{\varepsilon_i }{\|K\|_i}\right),\ \ \text{for all   } i \in \{1,...,n\}.$$ Then,
 $$S\cdot K\subseteq K'\cdot K+{B}\left(\theta,\mu(S)+r \|K^{-1}\|\right)\cdot K, $$
 where
 $$ r \|K ^{-1}\|:=\left(%
\begin{array}{ccc}
  \dfrac{r_1 }{\|K\|_1}  \\
  \vdots \\
   \dfrac{r_ n}{\|K\|}_n \\
\end{array}%
\right).$$
So
$$S\cdot K\subseteq K'\cdot K+{B}\left(\theta, \| K\|\mu(S)+r \right).$$
% from which, we infer that
 % $$S.K\subseteq K'.K+\overline{B_i}\big(0,\mu_i(S)+\dfrac{\epsilon_i }{\|K\|_i}\big)\|K\|_i, \ \ %\text{for all   } i \in \{1,...,n\}.$$
 Keeping in mind the subadditivity of
the  measure of weak noncompactness and using lemma \ref{lemma}, we get
   %Now, by , we have
   $$\begin{array}{rcl}\displaystyle\mu(S\cdot K)&\leq&\displaystyle\mu\left(K\cdot K'\right) +\mu\left({B}\left(0,\| K\|\mu(S)+r \right) \right)\\\\&\leq& \|K\|\mu(S)+r.\end{array}$$
  % $$\mu_i()\leq ,\ \ \text{for all   } i \in \{1,...,n\}.$$
    Since $\varepsilon$ is arbitrary, we deduce that $$\mu(S\cdot K)\leq \|K\|\mu(S).$$

\end{pf}

\bigskip

\begin{Lemma}\label{prop}Let $X$ be a generalized Banach algebra. Assume that
 $T:X \longrightarrow X$ is  weakly sequentially continuous.   If $T$ is contractive with Lipschitz matrix $M \in \mathcal{M}_{n \times n}(\mathbb{R}_+),$  then {for any bounded subset   } $S$ of $X,$ one has
  $$\mu\left(T(S)\right)\leq M\mu\left(S\right).$$
 % where,
 % \vskip0.2cm\hskip4cm$M=(a_{ij})_{1\leq i,j\leq n}.$
 %$A$ is $M$-set-contraction with respect to $\mu.$
 \hf
\end{Lemma}
\begin{pf}Let $S$ be a bounded subset of $X$  and  $\varepsilon \in \mathbb{R}^{n}$ such that $\varepsilon> \mu(S)$ and let $M=(a_{ij})_{1\leq i,j \leq n}.$  By the definition of $\mu,$ we have
%let $l_i\geq 0 $  such that $$\mu_i(S)= l_i,  \text{for all }\ \ i \in \{1,...,n\}. $$
there exists $r\in \mathbb{R}^{n}$   and a weakly compact subset $K$ of $X$ such that $0<r<\varepsilon$ and
$S \subseteq  K+{B}(\theta,{r}).$
Let $y\in T\left(K+{B}(\theta,{r})\right),$ then there exists $x\in K+{B}(\theta,{r})$ such that $y=Tx.$ Since $x\in K+{B}(\theta,{r}),$ there are  $k \in K$ and $b\in {B}\left(\theta,{r})\right)$ such that $x=k+b,$ and so
$$\|y-Tk\|_i=\|Tx-Tk\|_i\leq \sum_{j=1}^{n} a_{ij}\|x-k\|_j=\sum_{j=1}^{n} a_{ij}\|b\|_j\leq \sum_{j=1}^{n} a_{ij}r_j.$$
This means that
$$\|y-Tk\|\leq \left(
                 \begin{array}{c}
                   \displaystyle\sum_{j=1}^{n} a_{1j}r_j\\
                    \displaystyle\sum_{j=1}^{n} a_{2j}r_j\\
                   \vdots \\
                    \displaystyle\sum_{j=1}^{n} a_{nj}r_j\\
                 \end{array}
               \right)
=Mr.$$
That is, $y\in TK+{B}(\theta,{Mr})$ and consequently $TS\subset TK+{B}(\theta,{Mr}).$ Moreover, since $T$ is sequentially weakly continuous we have $\overline{TK}^{w}\in \mathcal{W}(X).$ Accordingly, $$\mu(TS)\leq Mr.$$
Letting $\varepsilon \to \mu(S)$ in the above inequality, we get  $\mu(TS)\leq M\mu(S).$
\end{pf}

\bigskip

\noindent \textbf{ Proof of Theorem \ref{3}.} Following the same procedures as in the proof of Theorem \ref{3.1}, it can be proved that the inverse operator
$T:=\left(\displaystyle\dfrac{I-C}{A}\right)^{-1}B$ exists on $S.$ Now, we claim that $T(S)$ is a relatively compact subset of $X.$ If this is not the case, then $\mu(TS)>0_{\mathbb{R}^{n}}.$ Keeping in mind the subadditivity of
the De Blasi's measure of weak noncompactness and using  the  equality:
\begin{eqnarray}\label{dec}
 %\nonumber to remove numbering (before each equation)
 \displaystyle
 T=AT\cdot B+CT
 \end{eqnarray}
we obtain
 % In order to acheive the proof, we will apply Theorem \ref{1}. Hence, we only have to prove that the operator $T: S\longrightarrow S$ is weakly sequentially continuous and $T(S)$ is relatively weakly compact. For this purpose,  taking into
%account the subadditivity of
%the  and using the following equality  to obtain
  %using both the equality  and knowing the subadditivity of the De Blasi’s measure of weak noncompactness, we get
\begin{align*}
\mu(T(S)) &\leq \mu \left(AT(S)\cdot B(S)\right)+\mu \left(CT(S)\right).
\end{align*}
The use of Lemma \ref{lemma3.2} as well as Lemma \ref{prop} leads to
%Now, by  and proposition \ref{prop},  we obtain
%the equality a
%\begin{equation}\label{eq1}
%T= ATB+CT,
%\end{equation}
\begin{align*}
\mu(TS)&\leq \left( \|B(S)\|M_A+M_C\right) \mu(TS).
\end{align*}
Since $\|B(S)\|M_A+M_C $ is a matrix converging to zero, we get a contradiction and consequently the claim is approved.
%From lemma \ref{lem},  we infer that $$\mu(T(S))\leq \big(I-bM_A-M_C \big)^{-1}0_{\mathbb{R}_{+}^{n}}.$$ So $$\mu(T(S))=0_{\mathbb{R}_{+}^{n}}.$$ Then, $T(S)$ is relatively weakly compact. Next, a standard argument $($see Theorem \ref{3.1} $)$ guarantees that $T$ is sequentially weakly continuous and the result is proved.\
\hfill{Q.E.D.}

\bigskip

\noindent As easy consequences of Theorem \ref{3} we obtain the following result.
\bigskip
\begin{Corollary}
 Let $S$ be a nonempty, bounded, closed and convex subset of  a generalized Banach algebra $X$ satisfying the sequential condition  $(\mathcal{GP}).$ Assume that $A,C:X\longrightarrow X$ and $B:S \longrightarrow X$  are three weakly sequentially continuous operators satisfying the following conditions:\\
\noindent $(i)$ $A$ is regular and is contractive with a Lipschitz matrix $M,$\\
\noindent $(ii)$ $B$ and $C$ are weakly compact, \\
\noindent $(iii)$ $\left(\displaystyle\dfrac{I-C}{A}\right)^{-1}$ exists on $B(S),$\\
\noindent $(iv)$ $x=Ax\cdot By+Cx,\ \ y \in S \Rightarrow  x \in S. $\\
 Then, the operator equation $($\ref{abc}$)$ has, at least, one solution in $S$ as soon as $\|B(S)\|M$ is a matrix converging to zero. \hf
\end{Corollary}
\bigskip
\noindent In the following result, we will consider that $A$ and $B$ are weakly compact operators.

\bigskip
\begin{Theorem}\label{rem} Let $S$ be a nonempty, bounded, closed and convex subset of  a generalized Banach algebra $X$ satisfying the condition  $(\mathcal{GP}).$ Assume that $A, B, C:S\longrightarrow X$ are three weakly sequentially continuous operators satisfying the following conditions:\\
\noindent{$(i)$} $A$ is regular,\\
\noindent{$(ii)$} $A$ and $B$ are weakly compact,\\
\noindent {$(iii)$} If $(I-C)x_n\rightharpoonup y,$ then there exists a weakly convergent subsequence of $(x_n)_n,$\\
\noindent{$(iv)$} $\left(\displaystyle\dfrac{I-C}{A}\right)^{-1}$ exists on $B(S)$, and\\
\noindent{$(v)$} $x=Ax\cdot By+Cx \in S,\ \ y\in S$  $\Rightarrow$ $x \in S.$\\
 Then, the operator equation $($\ref{abc}$)$ has, at least, one solution in $S.$ \hf
\end{Theorem}
\begin{pf}
%From assumption $(iii),$ we deduce that $T:=\left(\displaystyle\dfrac{I-C}{A}\right)^{-1}B: S\longrightarrow S $ is well defined.
% Now, we prove that $T(S)$ is relatively weakly compact.
 Let $\{y_n, \ n\in \mathbb{N} \}\subset T(S),$ there is a sequence  $\{x_n, \ n\in \mathbb{N} \}\subset S$ such that $$y_n=Tx_n=\left(\displaystyle\dfrac{I-C}{A}\right)^{-1}Bx_n.$$ Or equivalently, $y_n=Ay_n\cdot Bx_n+Cy_n.$ Taking into
account the weak compactness of  the weak closure of $A(S)$ and $B(S),$ we infer that
  %to obtain there is two subsequence   Since $A(S)$ and $B(S)$ are relatively weakly compact, it follows that there exist $(y_{n_{k}})_k$ of $(y_n)_n$ and $(x_{n_{k}})_k$ of $(x_n)_n$ such that
  $$ATx_{n_{k}}\rightharpoonup x  ~~~~~~~~~~~~~~~~\textrm{  and }~~~~~~~~~~~~~~~~ Bx_{n_{k_{j}}}\rightharpoonup x', \textrm{ for some } x,x'\in S,$$
 where $(x_{n_{k}})$ is a subsequence of  $\{x_n, \ n\in \mathbb{N} \}$ and $(x_{n_{k_{j}}})$ is a subsequence of  $\{x_{n_{k}}, \ k>n \}.$ Using the condition $(\mathcal{GP}),$ we get
  $$Ay_{n_{k_{j}}}\cdot Bx_{n_{k_{j}}}=(I-C)y_{n_{k_{j}}}\rightharpoonup x\cdot x' \ \ \text{ in } X.$$  Based on assumption $(iii),$ it follows that  there exists a weakly convergent subsequence of $(y_{n_{k_{j}}})$ and consequently $T(S)$ is a relatively weakly compact subset of $X.$  The use of Theorem \ref{3.1} achieves the proof. %  An argument similar to that in the proof of Theorem , leads to the weak sequential continuity of the operator $T.$ The result follows from Theorem \ref{1}.
\end{pf}

\bigskip

\begin{remark}If we assume that $C$ is sequentially weakly continuous and contractive with a Lipschitz Matrix $M,$  then $C$ satisfies the condition $(iii)$ of Theorem \ref{rem}. In fact, let $\{x_n, \ n\in \mathbb{N} \}$ be a sequence in $S$ such that  $(I-C)x_n\rightharpoonup y,$ for some  $y \in X.$ Based on the subadditivity of the De Blasi's measure of weak non-compactness it is shown that
$$\beta(\{x_n, \ n\in \mathbb{N}  \})\leq \beta(\{(I-C)x_n, \ n\in \mathbb{N} \})+\beta(\{Cx_n,  \ n\in \mathbb{N} \}). $$
If we consider the weak compactness of the weak closure of $\{(I-C)x_n, \ n\in \mathbb{N} \}$  and deploy Lemma  \ref{prop}, we get $$\beta(\{x_n, \ n\in \mathbb{N}  \})\leq M\beta(\{x_n, \ n\in \mathbb{N} \}).$$ Since $M$ is a matrix converging to zero, then there is a weakly convergent subsequence of $\{x_n, \ n\in \mathbb{N}  \}.$ \hf
%Taking into account that is relatively weakly compact and using proposition , we obtain
%$$\mu(\{x_n\})=0_{\mathbb{R}_{+}^{n}}. $$ This shows that, there exists a weakly convergent subsequence of $(x_n)_n.$
% Assume that $A,B,C:S\longrightarrow X$   are three weakly sequentially continuous operators on $S.$ Suppose that $A,B$ and $C$ satisfies conditions $(i),(ii), (iv)$ and $(v)$ of Theorem \ref{rem} and $C$ is contractive with a Lipschitz Matrix $M.$ Then the conclusion of Theorem \ref{rem} holds. \hf
%Let $S$ be a bounded subset of $X.$ Every $M$- contractive and weakly sequentially continuos mapping $C: S\longrightarrow X,$  satisfies the assumption $(iii)$ of Theorem \ref{rem}
% if $$(I-A)x_n \rightharpoonup y,$$ then there exists a weakly convergente subsequence of $(x_n)_n.$ \hf
\end{remark}

\bigskip
\noindent A consequence of Theorem \ref{rem} is
\bigskip
%\par We now consider the case when $A(S)$ is not relatively weakly compact. We obtain the following result
\begin{Corollary}Let $S$ be a nonempty, bounded, closed and convex subset of  a generalized Banach algebra $X$ satisfying the condition  $(\mathcal{GP}).$ Assume that $A,B,C:S\longrightarrow X$ are three weakly sequentially continuous operators satisfying the following conditions:\\
\noindent{$(i)$} $A$ is regular and $B$ is weakly compact,\\
\noindent {$(ii)$} If $\left(\displaystyle\dfrac{I-C}{A}\right) x_n\rightharpoonup y,$ then there exists a weakly convergent subsequence of $(x_n)_n,$\\
\noindent{$(iii)$} $\left(\displaystyle\dfrac{I-C}{A}\right)^{-1}$ exists and $\left(\displaystyle\dfrac{I-C}{A}\right)^{-1}B(S)\subseteq S.$\\
%\noindent{$(iv)$} $Ax{\cdot}Bx+Cx\in S,\ \ \text{for all } x \in S.$\\
 Then, the  operator equation $($\ref{abc}$)$ has, at least, one solution in $S.$ \hf
\end{Corollary}
\bigskip
\noindent Let us study the case where the operators $B$ and $\left(I-\displaystyle\dfrac{I-C}{A}\right)$ are $M_1$-$\mu$-contraction and  $M_2$-$\mu$-contraction respectively.
\bigskip
\begin{Theorem}\label{mesure} Let $S$ be a nonempty, bounded, closed and convex subset of  a generalized Banach algebra $X.$ Assume that $A, C:X\longrightarrow X$ and $B:S\longrightarrow X$ are three operators satisfying the following conditions:\\
\noindent{$(i)$} $A$ is regular,\\
\noindent{$(ii)$} $B$ and $\left(I-\displaystyle\dfrac{I-C}{A}\right)$ are $M_1$-$\mu$-contraction and $M_2$-$\mu$-contraction respectively, \\
 \noindent{$(iii)$} $B$ and $\left(\displaystyle\dfrac{I-C}{A}\right)$ are weakly sequentially continuous,\\
%\noindent {$(iv)$} $\left(I-\displaystyle\dfrac{I-C}{A}\right)$ is $M_2$-$\mu$-contraction,\\
\noindent{$(iv)$} $\left(\displaystyle\dfrac{I-C}{A}\right)^{-1}$ exists on $B(S),$ and\\
\noindent{$(v)$} $x=Ax\cdot By+Cx,\ \ y\in S$  $\Rightarrow$ $x \in S.$\\
 Then, the operator equation $($\ref{abc}$)$ has, at least, one solution in $S$ as soon as $(I-M_2)^{-1}M_1 $ is a matrix converging to zero.\hf
\end{Theorem}
\begin{pf}
It is easy to see that the operator $T:=\left(\displaystyle\dfrac{I-C}{A}\right)^{-1}B:S\longrightarrow S$ is well defined. %In order  to apply Theorem \ref{condensing}, it is sufficient to demonstrate that the operator $T$ is weakly sequentially continuous and $\mu$-condensing. To do so, consider
Let  $\{x_n,\ \ n \in \mathbb{N}\}$ be a weakly convergent sequence of $S$ to a point $x\in S.$ Keeping in mind the weak sequential continuity of the operator $B$ and using the equality
\begin{eqnarray}\label{mu}
T=B+\left(I- \dfrac{I-C}{A}\right)T,
\end{eqnarray}
 we obtain
$$\begin{array}{rcl}\mu \left(\{Tx_n, \ \ n \in \mathbb{N}\}\right)&\leq&\displaystyle\mu\left(Bx_n, \ \ n \in \mathbb{N}\right) +\mu\left( \left(I-\dfrac{I-C}{A}\right)(\{Tx_n,\ \ n \in \mathbb{N}\})\right)\\\\&\leq& M_2 \mu \left(\{Tx_n, \ \ n \in \mathbb{N}\}\right).
\end{array}$$
This inequality means, in particular, that $\{Tx_n, \ \ n \in \mathbb{N}\}$ is relatively weakly compact. Consequently, there is a subsequence $(x_{n_k})_k$ of $\{Tx_n, \ \ n \in \mathbb{N}\}$  such that $Tx_{n_k}\rightharpoonup y,$ fore some $y \in S.$ Making use of equality $(\ref{mu}),$ together with the assumptions on $B$ and $\left(\frac{I-C}{A}\right)$; enables us to have $y = Tx.$
 Now a standard argument shows that $Tx_n\rightharpoonup Tx.$
    Suppose the contrary, then there exists  a weakly neighborhood $V^{w}$ of $Tx$ and a subsequence $(x_{n_j})_j$  of  $\{x_{n}, \ n\in \mathbb{N} \}$  such that $Tx_{n_{j}} \notin V^{w},$ for all $j\geq 1.$  Arguing as
before, we may extract a subsequence  $(x_{n_{j_{k}}})_{k \in \mathbb{N}}$ of $\{x_{n_j}, \ j\in \mathbb{N} \}$ verifying  $Tx_{n_{j_{k}}}\rightharpoonup y,$  which is a contradiction and consequently $T$ is sequentially weakly continuous. Next, $T$ is $\mu$-condensing. In fact, let $V$ be a bounded subset of $S$ with $\mu(V) > 0.$  Using the subadditivity of the De Blasi's measure of weak
noncompactness, we get
$$\begin{array}{rcl}
\mu\left(T(V) \right) &\leq& \mu \left(B(V)\right)+ \mu \left(\left(I-\dfrac{I-C}{A}\right)T(V) \right) \\\\&\leq& M_1\mu(V)+M_2\mu(T(V)).
\end{array}$$
This implies that
$$\mu\left(T(V) \right)\leq (I-M_2)^{-1}M_1\mu\left(V\right). $$
Hence, $T$ has, at least, one fixed point $x$ in $S$ in view of Theorem \ref{condensing} .
\end{pf}

\bigskip

\noindent If we take $A = 1_X$ in the above result, where $1_X$ is the unit element of the generalized Banach algebra $X$, we obtain the following Corollary.

\bigskip
\begin{Corollary}
Let $S$ be a nonempty, bounded, closed and convex subset of  a generalized Banach algebra $X.$ Assume that $ C:X\longrightarrow X$ and $B:S\longrightarrow X$ are two sequentially weakly continuous operators satisfying the following conditions:\\
\noindent{$(i)$} $B$ and $C $ are  $M_1$-$\mu$-contraction and $M_2$-$\mu$-contraction respectively,\\
%\noindent {$(ii)$} $C $ is $M_2$-$\mu$-contraction,\\
\noindent{$(ii)$} $\left(I-C \right)^{-1}$ exists on $B(S),$\\
\noindent{$(iii)$} $x= By+Cx,\ \ y\in S$  $\Rightarrow$ $x \in S.$\\
 Then, the operator equation $($\ref{abc}$)$ has, at least, one solution in $S$ as soon as $(I-M_2)^{-1}M_1 $ is a matrix converging to zero.\hf
\end{Corollary}
\begin{remark}
Note that condition $(iii)$ in Theorem $\ref{mesure}$ may be replaced by "$A,B$ and $C$ are weakly sequentially continuous", but the generalized Banach algebra must satisfy condition $(\mathcal {GP}).$ Now, we can study the folllwing result. \hf
\end{remark}
\begin{Theorem} Let $S$ be a nonempty, bounded, closed, and convex subset of a generalized Banach algebra $X$ satisfying the sequential condition  $(\mathcal{GP}).$ Assume that  $A, C:X \longrightarrow X$ and $B:S \longrightarrow X$ are three weakly sequentially  continuous operators satisfying the following conditions:\\
\noindent{$(i)$}  $A$ is regular and weakly compact,\\
\noindent{$(ii)$} $B$ and $C $ are  $M_1$-$\mu$-contraction and $M_2$-$\mu$-contraction respectively, \\
%\noindent {$(iii)$} $C $ is $M_2$-$\mu$-contraction,\\
\noindent{$(iii)$} $\left(\displaystyle\dfrac{I-C}{A}\right)^{-1}$ exists on $B(S),$ and\\
\noindent{$(iv)$} $x=Ax\cdot By+Cx,\ \ y\in S$  $\Rightarrow$ $x \in S.$\\
Then, the operator equation $($\ref{abc}$)$ has, at least, one solution in $S$ as soon as $(I-M_2)^{-1}\|A(S)\|M_1 $ is a matrix converging to zero.\hf

\end{Theorem}

\section{Integral Equations of Fractional Order}
Let $(X,\|\cdot\|)$ be a reflexive Banach  and let ${C}(J, {X})$ be the Banach algebra of all
$X$-valued continuous functions defined on $J=[0,1]$, endowed with the norm $\|f\|_\infty = \displaystyle\sup_{t\in J}\|f(t)\|.$
We will use Theorem \ref{3} to examine the existence of solutions to the coupled system  of quadratic integral equations of fractional order $(\ref{Sy3}).$
%\begin{equation}
%x_i(t)= f_i(t,x_1(t),x_2(t))\cdot \displaystyle\int_{0}^{t}\displaystyle\frac{(t-s)^{\alpha_i-1}}%{\Gamma(\alpha_i)}g_i(s,x_1(s),x_2(s))ds+\displaystyle\sum _{k=1}^{m}I^{\beta_{i}^{k}}h_{i}^{k}%(t,x_1(t),x_2(t)),\ \ i=1,2$$
%\end{equation}
%where $J=[0,1],$ and $x_i=x_i(t)$ are unknown functions.
  We need the following definition and proposition  in the sequel.
\begin{Definition} \cite{Pettis}
Let $f:[0,T]\rightarrow X $ be a function. The fractional Pettis integral of the
function $f$ of order $ \alpha \in  \mathbb{R_+}$  is defined by
$$ I^{\alpha}f(t)=\int_{0}^{t}\frac{(t-s)^{\alpha -1}}{\Gamma(\alpha)}f(s)ds,  $$
where the sign $ "\int " $  denotes the Pettis integral.
\end{Definition}
\begin{Proposition}\cite{Pettis}
If $f:[0,T]\rightarrow X $ is Riemann integrable on $[0,T],$ then $I^{\alpha}f$ exists on $[0,T]$  and fractional Pettis integral.

\end{Proposition}

\noindent Let us now introduce the following assumptions:\\
\noindent{($\mathcal{H}_0$)} The function $f_i:J\times X \times X \longrightarrow X, i=1,2 $ is such that:\\
\hspace*{20pt} {$($a$)$} The partial function $x\longrightarrow f_1(t,x,y)$ is regular on $X$% i.e. ${f}^{1}_{(t,y(t))}(x(t))=\dfrac{x(t)}{f_1(t,x(t),y(t))}$ is well defined.
\\
\hspace*{20pt} {$($b$)$} The partial function $y\longrightarrow f_2(t,x,y)$ is regular on $X$% i.e. ${f}^{2}_{(t,x(t))}(y(t))=\dfrac{y(t)}{f_2(t,x(t),y(t))}$ is well defined .
\\
%The functions $(x,y)\rightarrow \dfrac{x}{f_1(t,x,y)}$ and $(x,y)\rightarrow \dfrac{y}{f_2(t,x,y)}$ are well defined.\\
\hspace*{20pt} {$($c$)$} The partial function $t \mapsto f_i(t,x,y) $ is continuous,\\
\hspace*{20pt} {$($d$)$}  The partial function $(x,y) \mapsto f_i(t,x,y)$ is weakly sequentially continuous,\\
\hspace*{20pt} {$($e$)$} There is  nonnegative real numbers  $ a_{i1}$ and $a_{i2}  ,$  $i=1,2$  such that \\
$$\|f_i(t,x,y)-f_i(t,\tilde {x},\tilde{y})\|\leq a_{i1}\|x-\tilde{x}\|+a_{i2}\|y-\tilde{y}\|.$$
\noindent{($\mathcal{H}_1$)} The function $g_i:J\times X \times X \longrightarrow X, i=1,2 $ is such that:\\
\hspace*{20pt} {$($a$)$}The partial function $t \mapsto g_i(t,x,y)$ is Riemann integrable,\\
\hspace*{20pt} {$($b$)$} The partial function $(x,y) \mapsto g_i(t,x,y)$ is weakly sequentially  continuous.\\
%\hspace*{20pt} {$($c$)$}There exists $P > 0$ such that for all $r_0 > 0;$
% $$$$

%There exists a function $P \in C(J,\mathbb{R}_{+})$ and there is  two continuous nondecreasing functions  $\psi_1,\psi_2:[0,\infty)\longrightarrow [0,\infty)$ such that:
%$$\|g_i(t,x(t),y(t))\|\leq P(t)[\psi_1(\|x\|)+\psi_2(\|y\|)].$$
\noindent{($\mathcal{H}_2$)} The function $h_{i}^{k}:J \times X \times X \longrightarrow X,$  $k= {1,\ldots,m}$ is such that:\\
\hspace*{20pt} {$($a$)$}The partial function $t \mapsto  h_i^k(t,x,y)$ is Riemann integrable,\\
\hspace*{20pt} {$($b$)$} The partial function $(x,y) \mapsto h_i^k(t,x,y)$ is weakly sequentially  continuous,\\
\hspace*{20pt} {$($c$)$}  There is nonnegative real numbers  ${b}_{i1}^{k}$ and ${b}_{i2}^{k}$  such that \\
$$\|h_i^k(t,x,y)-h_i^k(t,\tilde {x},\tilde{y})\|\leq {b}_{i1}^{k}\|x-\tilde{x}\|+{b}_{i2}^{k}\|y-\tilde{y}\|.$$
%\noindent{($\mathcal{H}_3$)} There exists a number $r>0$ such that:\\
%
%\begin{equation} \label{r}
%r\geq \frac{F_0\|P\|_\infty \big(\psi_1(r)+\psi_2(r)\big)\big(\frac{T^{\alpha_1}}{\Gamma(\alpha_1+1)}+\frac{T^{\alpha_2}}{\Gamma(\alpha_2+1)}\big)+H_0\sum_{k=1}^{m}\big(\frac{T^{\beta_1^k}}{\Gamma(\beta_1^k+1)}+\frac{T^{\beta_2^k}}{\Gamma(\beta_2^k+1)}\big)
% }{1-C \big[\|P\|_\infty \big(\psi_1(r)+\psi_2(r) \big) \big(\frac{T^{\alpha_1}}{\Gamma(\alpha_1+1)}+\frac{T^{\alpha_2}}{\Gamma(\alpha_2+1)} \big)+\sum_{k=1}^{m} \big(\frac{T^{\beta_1^k}}{\Gamma(\beta_1^k+1)}+\frac{T^{\beta_2^k}}{\Gamma(\beta_2^k+1)}\big)\big]} ,
%\end{equation}
% where $$F_0=\max\ \big\{ \|f_1(\cdot,0,0)\|,\|f_2(\cdot,0,0)\| \big\},$$
% $$ H_0= \max\ \big\{ \|h_1^k(\cdot,0,0)\|,\|h_2^k(\cdot,0,0)\|;\ \ k=1,...,m \big\},$$
%        % $$ H_0= \sup\limits_{ 0\leq t\leq T} \|h_j^i(t,0,0)\|, \ \ i=1,...,n\ \ \ ;j=1,2,$$
%$$C=\max\{ a_{11},a_{12},a_{21},{a_{22}},{b}_{11}^{k},{b}_{12}^k,{b}_{21}^K,{b}_{21}^k,{b}_{22}^k ;\ \ k=1,...,m  \}, $$
%and
%$$ C\big[\|P\|_\infty \big(\psi_1(r)+\psi_2(r) \big) \big(\frac{T^{\alpha_1}}{\Gamma(\alpha_1+1)}+\frac{T^{\alpha_2}}{\Gamma(\alpha_2+1)} \big)+\sum_{k=1}^{m} \big(\frac{T^{\beta_1^k}}{\Gamma(\beta_1^k+1)}+\frac{T^{\beta_2^k}}{\Gamma(\beta_2^k+1)}\big)\big]< 1. $$\\
\begin{Theorem} \label{mple}  Suppose that the assumptions $(\mathcal{H}_{0})-(\mathcal{H}_3)$ are satisfied. Moreover, assume that there exists a real number $r_0> 0$ and $P \in \mathbb{R}^{\ast}_+$ such that
\begin{eqnarray}\label{app}
 \displaystyle
\left\{
    \begin{array}{lll}
     \|g_i(t, s, x, y)\|\leq P  \ \ \text{   with  } \|x\|\leq r_0 \text{  and   }\|y\|\leq r_0 \\\\
      \rho\big[ P\big(\frac{T^{\alpha_1}}{\Gamma(\alpha_1+1)}+\frac{T^{\alpha_2}}{\Gamma(\alpha_2+1)} \big)+\sum_{k=1}^{m} \big(\frac{T^{\beta_1^k}}{\Gamma(\beta_1^k+1)}+\frac{T^{\beta_2^k}}{\Gamma(\beta_2^k+1)}\big)\big]< 1, \textrm{ and } \\\\
M_A, M_{C} \textrm{ and } \|B(S)\|M_A+M_C \textrm{ are three matrices converging to zero, where }
    \end{array}
  \right.
\end{eqnarray}
%Moreover,
$$\rho=\max\{ a_{11},a_{12},a_{21},{a_{22}},{b}_{11}^{k},{b}_{12}^k,{b}_{21}^K,{b}_{21}^k,{b}_{22}^k ;\ \ k=1,\ldots,m  \},$$
$$M_A=\left(
\begin{array}{cc}
  a_{11}  & a_{12} \\
  a_{21} & a_{22}   \\
\end{array}
\right) , \ \ M_{C}=\left(
\begin{array}{cc}
 \sum_{k=1}^{m}\frac{T^{\beta_1^{k}}}{\Gamma(\beta_1^{k}+1)}b{_{11}^{k} } & \sum_{k=1}^{m}\frac{T^{\beta_1^{k}}}{\Gamma(\beta_1^{k}+1)}b{_{12}^{k} } \\
  \sum_{k=1}^{m}\frac{T^{\beta_2^{k}}}{\Gamma(\beta_2^{k}+1)}b{_{21}^{k} } & \sum_{k=1}^{m}\frac{T^{\beta_2^{k}}}{\Gamma(\beta_2^{k}+1)}b{_{22}^{k} }   \\
\end{array}
\right) \texttt{ and }$$ $$\ \ \|B(S)\|M_A+M_C= \left(
                                   \begin{array}{cc}
                                     \dfrac{PT^{\alpha_1}}{\Gamma(\alpha_1+1)}a_{11}+\sum_{k=1}^{m}\frac{T^{\beta_1^{k}}}{\Gamma(\beta_1^{k}+1)}b{_{11}^{k} } & \dfrac{PT^{\alpha_2}}{\Gamma(\alpha_2+1)}a_{12}+\sum_{k=1}^{m}\frac{T^{\beta_2^{k}}}{\Gamma(\beta_2^{k}+1)}b{_{12}^{k} } \\
                                     \dfrac{PT^{\alpha_1}}{\Gamma(\alpha_1+1)}a_{21}+\sum_{k=1}^{m}\frac{T^{\beta_1^{k}}}{\Gamma(\beta_1^{k}+1)}b{_{21}^{k} } & \dfrac{PT^{\alpha_2}}{\Gamma(\alpha_2+1)}a_{22}+\sum_{k=1}^{m}\frac{T^{\beta_2^{k}}}{\Gamma(\beta_2^{k}+1)}b{_{22}^{k}} \\
                                   \end{array}
                                 \right).$$
%\left(
%\begin{array}{cc}
%  m_1  & m_{2} \\
%  m_{3} & m_{4}   \\
%\end{array}
%\right)  \in \mathcal{M}_{(2 \times 2)}(\mathbb{R}_+)$ are three matrices converging to zero, where $$m_1= \dfrac{T^{\alpha_1}P}{\Gamma(\alpha_1+1)}a_{11}+\sum_{k=1}^{2}\frac{T^{\beta_1^{k}}}{\Gamma(\beta_1^{k}+1)}b{_{11}^{k} } ,$$
%$$m_2= \dfrac{T^{\alpha_2}P}{\Gamma(\alpha_2+1)}a_{12}+\sum_{k=1}^{2}\frac{T^{\beta_2^{k}}}{\Gamma(\beta_2^{k}+1)}b{_{12}^{k} },$$
%$$m_3= \dfrac{T^{\alpha_1}P}{\Gamma(\alpha_1+1)}a_{21}+\sum_{k=1}^{2}\frac{T^{\beta_1^{k}}}{\Gamma(\beta_1^{k}+1)}b{_{21}^{k} },$$
%$$m_4= \dfrac{T^{\alpha_2}P}{\Gamma(\alpha_2+1)}a_{22}+\sum_{k=1}^{2}\frac{T^{\beta_2^{k}}}{\Gamma(\beta_2^{k}+1)}b{_{22}^{k}}.$$
Then the problem
$(\ref{Sy3})$ has a solution in $C(J,X)\times C(J,X).$ \hf

\end{Theorem}

\begin{pf}
We recall that the problem $(\ref{Sy3} $) is equivalent to the operator equation
$$(x,y)=(A_1(x,y),A_2(x,y))(B_1(x,y),B_2(x,y))+(C_1(x,y),C_2(x,y)),$$ where the operators $A_i,B_i$ and $C_i, \ \ i=1,2$ are defined by
$$\left\{
  \begin{array}{ll}
    A_i(x,y)(t)=f_i(t,x(t),y(t))\\\\
     B_i(x,y)(t)= \int_{0}^{t}\frac{(t-s)^{\alpha_1-1}}{\Gamma(\alpha_1)}g_i(s,x(s),y(s))ds\\\\
     C_i(x,y)(t)=\sum _{k=1}^{m}I^{\beta_{1}^{i}}h_{i}^{k}(t,x(t),y(t)).
  \end{array}
\right.$$

\noindent Let us define the subset $S$ of $C(J,X)\times C(J,X)$ by
 $$S=\left\{(x,y)\in C(J,X)\times C(J,X), \ \ \|(x,y)\| = \left(
\begin{array}{cc}
  \|x\|_\infty  \\
   \|y\|_\infty \\
\end{array} \right) \leq \left(%
\begin{array}{cc}
  r_0  \\
  r_0\\
\end{array}%
\right) \right\},$$ where
$$
r_0\geq\frac{F_0P\big(\frac{T^{\alpha_1}}{\Gamma(\alpha_1+1)}+\frac{T^{\alpha_2}}{\Gamma(\alpha_2+1)}\big)+H_0\sum_{k=1}^{m}\big(\frac{T^{\beta_1^k}}{\Gamma(\beta_1^k+1)}+\frac{T^{\beta_2^k}}{\Gamma(\beta_2^k+1)}\big)
 }{1-\rho\big[P \big(\frac{T^{\alpha_1}}{\Gamma(\alpha_1+1)}+\frac{T^{\alpha_2}}{\Gamma(\alpha_2+1)} \big)+\sum_{k=1}^{m} \big(\frac{T^{\beta_1^k}}{\Gamma(\beta_1^k+1)}+\frac{T^{\beta_2^k}}{\Gamma(\beta_2^k+1)}\big)\big]} ,
$$
 with
$$\left\{
    \begin{array}{ll}
      F_0=\max\ \big\{ \|f_1(\cdot,0,0)\|,\|f_2(\cdot,0,0)\| \big\} \\
      H_0= \max\ \big\{ \|h_1^k(\cdot,0,0)\|,\|h_2^k(\cdot,0,0)\|;\ \ k=1,\ldots,m \big\}\\
    \end{array}
  \right.
$$
Our strategy is to apply Theorem \ref{3} to prove the existence of a fixed point for the nonlinear equation $(\ref{Sy3})$ in $S.$ Then, we need to verify the following steps:\\
\textbf{Claim 1:} We start by showing that the operators $A,C: C(J,X) \times C(J,X)\longrightarrow C(J,X) \times C(J,X)$ and $B: S\longrightarrow C(J,X) \times C(J,X)$ are weakly sequentially continuous. Firstly, we verify that the operator $A_i(x,y), \text{ for } i=1,2$  is continuous on $J$ for all  $(x,y) \in C(J,X) \times C(J,X).$ To see this, let  $\{t_n, \ \ n \in \mathbb{N}\}$ be any
sequence in J converging to a point $t$ in $J.$    Then,
$$\begin{array}{rcl}
\|A_i(x,y)(t_n)-A_i(x,y)(t)\| &=& \|f_i(t_n,x(t_n),y(t_n))-f_i(t,x(t),y(t))\|
\\\\&\leq &  a_{i1}\|x(t_n)-x(t)\|+a_{i2}\|y(t_n)-y(t)\| \\\\&+& \|f_i(t_n,x(t),y(t)-f_i(t,x(t),y(t)\|.
\end{array}$$
The continuity of $x,y$ and $t\mapsto f_i(t,x,y)$ on $[0,1]$
 implies that the function $A_i(x,y)$ is continuous.
% In the proof of Theorem \ref{3}, we need to prove that $A$ and $C$ are weakly sequentially continuous on $E\times E$ and B is weakly sequentially continuous on $S.$ We begin to show the property for the operator $A.$ To see this,
Let $\{(x_n,y_n), n\in \mathbb{N}\}$ be a weakly convergent sequence of $C(J,X) \times C(J,X)$ to a point $(x,y).$ In this case, the set $\{(x_n,y_n), n\in \mathbb{N}\}$
 is bounded and so, we can apply the Dobrakov's theorem \cite {Dob} in order to get $$(x_n(t),y_n(t))\rightharpoonup (x(t),y(t)) \text{   in  } X\times X.$$  Based on assumption $(\mathcal{H}_0)(d),$ it is shown that $A_i(x_n,y_n)(t)\rightharpoonup A_i(x,y)(t)$ and  then, we can again apply the Dobrakov's theorem to obtain the weak sequential continuity of the operator $A.$ Besides, the use of assumption  $(\mathcal{H}_{1})(b) $   and  assumption $(\mathcal{H}_2)(b)$ as well as the Dobrakov's theorem \cite[page 36]{Dob} leads to the two maps $B$ and $C$  are weakly sequentially continuous.\\
 %combined with  the dominated convergence theorem for Pettis integral, we obtain
%$$I^{\beta_i^k}h^k_i(t,x_n(t),y_n(t))\rightharpoonup I^{\beta_i^k}h^k_i(t,x(t),y(t))  \text{  in }  X,  \ \  \text{ for } i   \in \{1,2\} \text{  and  } k \in  \{1,...,m\}. $$
% Thus $C_i(x_n,y_n)(t)\rightharpoonup C_i(x,y)(t), $ for $i\in \{1,2\}.$
%  Using again Dobrakov's theorem \cite{Dob}, we have $ C_i(x_n,y_n)\rightharpoonup C_i(x,y),$ for $i\in \{1,2\}.$
% Hence, $C$ is weakly sequentially continuous.
%Moreover, let $\{(x_n,y_n) \}_n$ be any sequence in $S$ weakly converging to a point $(x,y).$ We can apply  to obtain $$(x_n(t),y_n(t))\rightharpoonup (x(t),y(t))\ \ \text{   in   } X\times X.$$ Using both the Lebesgue dominated convergence theorem for Pettis integral and the assumption  we can deduce that $\ \ B_i(t,x_n(t),y_n(t))\rightharpoonup  B_i(t,x(t),y(t))\ \  \text{  in  } X,  i \in \{1,2\} .$
% Taking into account that
%$\{B_i(x_n,y_n), \ \ n \in \mathbb{N}  \}$ is bounded with bound $\dfrac{T^{\alpha_i}P}{\Gamma(\alpha_i+1)}$ and using the Dobrakov’s theorem  we show that $B$ is a weakly sequentially continuous operator on $S.$
\textbf{Claim 2:} The operators $A$ and $C$ are contractive. The claim regarding the operator $A$  is
immediate, from assumption $(\mathcal{H}_0)(e).$ Let us fix arbitrary $(x,y),(\tilde{x},\tilde{y}) \in C(J,X) \times C(J,X).$ If we
take an arbitrary $t\in[0, 1],$ then we get
$$\begin{array}{rcl} \|C_i(x,y)(t)-C_i(\tilde{x},\tilde{y})(t)\|&=& \left\|\displaystyle\sum_{k=1}^{m}I^{\beta_{i}^{k}}h_i^k(t,x(t),y(t)-\displaystyle\sum_{k=1}^{m}I^{\beta_{i}^{k}}h_i^k(t,\tilde{x}(t),\tilde{y}(t))\right\|  \\\\&\leq& \displaystyle\sum_{k=1}^{m}I^{\beta_{i}^{k}}\|h_i^k(t,x(t),y(t)-h_i^k(t,\tilde{x}(t),\tilde{y}(t))\| \\\\&\leq&  \displaystyle\sum_{k=1}^{m}I^{\beta_{i}^{k}}\left(b_{i1}^{k}\|x(t)-\tilde{x}(t)\| + b_{i1}^{k}\|y(t)-\tilde{y}(t)\| \right) \\\\&\leq& \displaystyle\sum_{k=1}^{m}\frac{T^{\beta_{i}^{k}}}{\Gamma(\beta_{i}^{k}+1)}\left(b_{i1}^{k}\|x(t)-\tilde{x}(t)\| + b_{i1}^{k}\|y(t)-\tilde{y}(t)\|  \right).   \end{array}$$
%By taking the Supremum over $t,$  we get
%$$ \|C_i(x,y)-C_i(\tilde{x},\tilde{y})\|_\infty \leq \displaystyle\sum_{k=1}^{m}\frac{T^{\beta_{i}^{k}}}{\Gamma(\beta_{i}^{k}+1)}\big(b_{i1}^{k}\|x-\tilde{x}\|_\infty + b_{i1}^{k}\|y-\tilde{y}\|_\infty  \big) ,\text{  for } i \in \{1,2\}.$$
This implies that
 $\left\|C(x,y)-C(\tilde{x},\tilde{y})\right\| \leq M_C\left\|(x,y)-(\tilde{x},\tilde{y})\right\|,$ where $$M_C= \left(
\begin{array}{cc}
 \displaystyle\sum_{k=1}^{m}\frac{T^{\beta_1^{k}}}{\Gamma(\beta_1^{k}+1)}b{_{11}^{k} } & \displaystyle\sum_{k=1}^{m}\frac{T^{\beta_1^{k}}}{\Gamma(\beta_1^{k}+1)}b{_{12}^{k} } \\ \\
  \displaystyle\sum_{k=1}^{m}\frac{T^{\beta_2^{k}}}{\Gamma(\beta_2^{k}+1)}b{_{21}^{k} } & \displaystyle\sum_{k=1}^{m}\frac{T^{\beta_2^{k}}}{\Gamma(\beta_2^{k}+1)}b{_{22}^{k} }
\end{array}
\right).$$ %This shows that $C$ is contractive with Lipschitz matrix $M_C $ and it is similarly clear that $A$ is contractive with Lipschtz matrix $M_A=\left(

\noindent \textbf{Claim 3:}  %$B$ is weakly compact on $S.$ Firstly, we begin to show that $B$ maps $S$ into $E \times E.$  If we take
Let $\varepsilon > 0,\ \ (x,y) \in S$ and $t,t' \in J$ such that
% $\big($without loss of generality assume that  $ \big(B_i(x,y)(t)-B_i(x,y)(t')\neq 0 \big)$ and
$|t'-t|<\varepsilon.$ From the Hahn-Banach theorem there exists a linear function $\phi \in X^{\ast}$ with $\|\phi\|=1$  such that $$\left\|B_i(x,y)(t')-B_i(x,y)(t)\right\|=\phi\left(B_i(x,y)(t)-B_i(x,y)(t)\right), \ \ i=1,2.$$ Based on the first inequality in \eqref{app}  it is shown that
$$\begin{array}{rcl}
\left\|B_i(x,y)(t')-B_i(x,y)(t)\right\| &=& \phi \left(\displaystyle\int_{0}^{t'}\frac{(t-s)^{\alpha_i-1}}{\Gamma(\alpha_j)}g_i(s,x(s),y(s))ds-\displaystyle\int_{0}^{t}\frac{(t-s)^{\alpha_i-1}}{\Gamma(\alpha_i)}g_i(s,x(s),y(s))\right) \\\\&\leq& \displaystyle\int_{0}^{t}\frac{|(t'-s)^{\alpha_i-1}-(t-s)^{\alpha_i-1}|}{\Gamma(\alpha_i)}\phi(g_i(s,x(s),y(s)))ds \\\\& +&  \displaystyle\int_{t}^{t'}\frac{(t'-s)^{\alpha_i-1}}{\Gamma(\alpha_i)}\phi(g_i(s,x(s),y(s)))ds  \\\\&\leq& P\displaystyle\int_{0}^{t}\frac{|(t'-s)^{\alpha_i-1}-(t-s)^{\alpha_i-1}|}{\Gamma(\alpha_i)} ds +  P\displaystyle\int_{t}^{t'}\frac{(t'-s)^{\alpha_i-1}}{\Gamma(\alpha_i)}ds.
\end{array}$$ This implies that
$\left\|B_i(x,y)(t')-B_i(x,y)(t)\right\|\to 0 \text{   as   } \varepsilon \to 0,$ and consequently
$B(S)$ is a  weakly equi-continuous subset.
%Thus, by the above findings, $$B(S)=\left\{\big(B_1(x,y),B_2(x,y)\big),  \ \ \|(x,y)\|_{E \times E}\leq \left(%
%\begin{array}{cc}
%  r_0  \\
%  r_0 \\
%\end{array}%
%\right)  \right\}$$ is continuous.
Let now $\{(x_n,y_n), n\in \mathbb{N}\}$ be any sequence in $S.$ From the first inequality in \eqref{app},  it follows that
$$\begin{array}{rcl}
\left\|B_i(t,x_n(t),y_n(t))\right\| &\leq& \displaystyle\int_{0}^{t}\frac{(t-s)^{\alpha_i-1}}{\Gamma(\alpha_i)}\left\|g_i(s,x_n(s),y_n(s))\right\|ds \\\\&\leq& \displaystyle\frac{PT^{\alpha_i}}{\Gamma(\alpha_i+1)},
\end{array}$$
for all $t\in [0,1].$
This demonstrate that $\{B_i(x_n,y_n), \ \ n \in \mathbb{N}\}$ is a uniformly bounded sequence in $B(S)$ and so,  $B(S)(t)$ is sequentially relatively weakly compact.  Hence, $B(S)$ is sequentially relatively
weakly compact in light of the Arzel\`{a}-Ascoli's theorem  \cite{Ascoli}. An application of Eberlein-\u{S}mulian's theorem \cite{kes}  yields
that $B(S)$ is relatively weakly compact.

 \noindent \textbf{Claim 4:} The operators $A,B$ and $ C $ satisfy assumption $(iv)$ of Theorem $\ref{3}.$ To see this, let $(x,y)\in C(J,X) \times C(J,X)$ and $(u,v) \in S$ with $(x,y)=A(x,y)\cdot B(u,v)+C(x,y).$ We shall show that $(x,y)\in S.$ For all $t\in[0,1],$ we have
%\left(a_{11}\|x(t)\|+a_{12}\|y(t)\|+\|f_1(t,0,0)\|\right)\frac{T^{\alpha_1}}{\Gamma(\alpha_1+1)}|P(t)|\left(\psi_1(r)+\psi_2(r)\right)\\\\&+&\displaystyle\sum_{k=1}^{m}\frac{T^{\beta_1^{k}}}{\Gamma(\beta_1^{k}+1)}\left(b_{11}^{k}\|x(t)\|+b_{12}^{k}\|y(t)\|+h_1^{k}(t,0,0) %\right) \\\\&\leq&
$$ \begin{array}{rcl}
\|x(t)\| &\leq& \left\|f_1(t,x(t),y(t))\right\|\displaystyle\int_0^t \displaystyle\frac{(t-s)^{\alpha_1-1}}{\Gamma(\alpha_1)}\left\|g_1(s,u(s),v(s))\right\|ds+ \displaystyle\sum_{i=1}^m I^{\beta_1^i}\|h_1^i(t,x(t),y(t))\| \\\\&\leq& \left[\rho( \|x\|_\infty+\|y\|_\infty)+F_0\right]\displaystyle\frac{PT^{\alpha_1}}{\Gamma(\alpha_1+1)}+\displaystyle\sum_{k=1}^{m}\frac{T^{\beta_1^{k}}}{\Gamma(\beta_1^{k}+1)}\left[\rho(\|x\|_{\infty}+\|y\|_{\infty})+H_0 \right]
\end{array}$$
This implies that
$$\begin{array}{rcl}
\|x(t)\|+\|y(t)\| &\leq& \big[\rho( \|x\|_\infty+\|y\|_\infty)+F_0\big]P\left(\frac{T^{\alpha_1}}{\Gamma(\alpha_1+1)}+ \frac{T^{\alpha_2}}{\Gamma(\alpha_2+1)}\right)\ \  \\  &+& \sum_{k=1}^{m}\left(\frac{T^{\beta_1^{k}}}{\Gamma(\beta_1^{k}+1)}+\frac{T^{\beta_2^{k}}}{\Gamma(\beta_2^{k}+1)}\right) \big[\rho(\|x\|_\infty+\|y\|_\infty)+H_0 \big]
\end{array}$$
Consequently, $$\|x\|_\infty\leq r_0 \text{   and   } \|y\|_\infty \leq r_0.$$
To end the proof, we apply Theorem $(\ref{mple})$, we deduce that  the problem $(\ref{Sy3})$ has,       at least,  one solution in $ C(J,X) \times C(J,X).$\\
\textbf{Example :}  Let $C(J,\mathbb{R})$ be the Banach algebra of all continuous functions from $J$ to $\mathbb{R} $ endowed with the sup-norm $\|\cdot\|$ defined by $\|x\|_\infty = \sup \limits_{0 \leq t \leq T}|x(t)|, $ for each $x \in C(J, \mathbb{R}).$

\begin{equation}\label{exe}
\left\{
  \begin{array}{lll}
  D^{\frac{1}{2}}\big[\frac{x(t)- \sum _{k=1}^{2}I^{\beta_{1}^{k}}h_{1}^{k}(t,x(t),y(t))}{f_1(t,x(t),y(t))}\big]=\frac{3}{35(13-t^2)}(7|x(t)|+15|y(t)|), & \hbox{} \\ \\

 D^{\frac{1}{2}}\big [\frac{y(t)- \sum _{k=1}^{2}I^{\beta_{2}^{k}}h_{2}^{k}(t,x(t),y(t))}{f_2(t,x(t),y(t))}\big]=\frac{3}{35(13-t^2)}(7|x(t)|+15|y(t)|), & \hbox{} \\ \\
 x(0)=y(0)=0, \ \ \hfill{ t  \in[0,1]}.
   \end{array}
\right.
\end{equation}
Note that this problem may be transformed into the fixed point problem $(\ref{Sy3})$ in view of lemma $2.5$ in \cite{D},
 where $$\sum_{k=1}^{2}I^{\beta_{1}^{k}}h_{1}^{k}(t,x(t),y(t))=I^{ 1 / 3}\frac{2te^{-3t}}{15(3+t)}\left(\frac{x^2(t)+9|x(t)|}{|x(t)|+5}+\frac{12e^{3t}}{5} \right) + I^{ 10/ 3}\frac{2t \sin \pi t}{14+t^2}\left(\frac{x^2(t)+5|x(t)|}{|x(t)|+8}+\frac{1}{3}\right),  $$
 $$\sum_{k=1}^{2}I^{\beta_{2}^{k}}h_{2}^{k}(t,x(t),y(t))=I^{ 7 / 4}\frac{t \sin t}{7(4+e^t)}\left(\frac{y^2(t)+4|y(t)|}{|y(t)|+3}+\cos t \right) + I^{ 29/ 6}\frac{3t \cos t}{10(4-t^2)}\left(\frac{y^2(t)+5|y(t)|}{|y(t)|+4}+\frac{t}{t+2}\right), $$
  and $$f_1(t,x(t),y(t))=\dfrac{3cos \pi t +2t}{5(2+10t^2)(|x(t)|+3)},f_2(t,x(t),y(t))=\dfrac{4cos \pi t +3t}{7(3+8t^2)(|y(t)|+6)},$$ and
  $$g_1(t,x(t),y(t))=g_2(t,x(t),y(t))=\frac{3}{35(13-t^2)}(7|x(t)|+15|y(t)|,$$
here $\alpha_1=\alpha_2 = \dfrac{1}{2}, T = 1, m = 2,  \beta_1^1=\dfrac{1}{3} ,\beta_1^2=\dfrac{10}{3}, \beta_2^1=\dfrac{7}{4}$ and $ \beta_2^2=\dfrac{29}{6}.$ \\We can show that
\begin{eqnarray*}
|f_1(t,x,y)-f_1(t,\tilde{x},\tilde{y})|\leq \frac{3+2t}{5(2+10t^2)}|x-\tilde{x}|,
\end{eqnarray*}
\begin{eqnarray*}
|f_2(t,x,y)-f_2(t,\tilde{x},\tilde{y})|\leq \frac{4+3t}{7(1+5t^2)}|y-\tilde{y}|,
\end{eqnarray*}
and
\begin{eqnarray*}
|h_1^1(t,x,y)-h_1^1(t,\tilde{x},\tilde{y})|\leq \frac{18 t}{75(3+t)}|x-\tilde{x}|,
\end{eqnarray*}
\begin{eqnarray*}
|h_1^2(t,x,y)-h_1^2(t,\tilde{x},\tilde{y})|\leq \frac{10 t}{8(14+t^2)}|x-\tilde{x}|,
 \end{eqnarray*}
 \begin{eqnarray*}
|h_2^1(t,x,y)-h_2^1(t,\tilde{x},\tilde{y})|\leq \frac{4 t}{21(4+e^t)}|y-\tilde{y}|,
\end{eqnarray*}
\begin{eqnarray*}
|h_2^2(t,x,y)-h_2^2(t,x,\tilde{y})|\leq \frac{3 t}{8(4-t^2)}|y-\tilde{y}|.
\end{eqnarray*}
It follows that $a_{11}=\dfrac{1}{12},\ \ a_{21}=0,\ \  a_{12}=0, \ \  a_{22}= \dfrac{ 1}{6},\ \ b^{1}_{11}=\dfrac{3}{50},\ \  b^{1}_{12}=0,\ \  b_{21}^{1}=0,\ \  b_{22}^{1}=\dfrac{4}{21(4+e)},\ \  b_{11}^{2}=\dfrac{1}{12},\ \ b_{12}^{2}=0,\ \  b_{21}^{2}=0,\ \  b_{22}^{2}= \dfrac{1}{8}. $
Now  $$|g_1(t,x,y)|=|g_2(t,x,y)|\leq \dfrac{1}{4} $$
%\dfrac{3 }{(13-t^2)}\left(\dfrac{|x|}{5}+\dfrac{3|y|}{7}\right),
%where $p(t)=\dfrac{3 }{(13-t^2)} $ \ \ and $\psi_1(x)=\dfrac{|x|}{5},\ \  \psi_2(y)=\dfrac{3|y|}{7}.$\\
It is easy to verify that $  \rho= \dfrac{1}{6},\ \  F_0=\dfrac{1}{36},\ \ H_0=\dfrac{2}{25}.$ We see that condition $(\mathcal{H}_3)$ is followed with a number $r_0=2.$      Moreover, $ \dfrac{PT^{\alpha_1}}{\Gamma(\alpha_1+1)}a_{11}+\sum_{k=1}^{m}\frac{T^{\beta_1^{k}}}{\Gamma(\beta_1^{k}+1)}b{_{11}^{k} }\simeq 0.0990,\ \
 \dfrac{PT^{\alpha_2}}{\Gamma(\alpha_2+1)}a_{12}+\sum_{k=1}^{m}\frac{T^{\beta_2^{k}}}{\Gamma(\beta_2^{k}+1)}b{_{12}^{k} }= \dfrac{PT^{\alpha_1}}{\Gamma(\alpha_1+1)}a_{21}+\sum_{k=1}^{m}\frac{T^{\beta_1^{k}}}{\Gamma(\beta_1^{k}+1)}b{_{21}^{k} }=0, \ \
  \dfrac{PT^{\alpha_2}}{\Gamma(\alpha_2+1)}a_{22}+\sum_{k=1}^{m}\frac{T^{\beta_2^{k}}}{\Gamma(\beta_2^{k}+1)}b{_{22}^{k}} \simeq 0.0653.$
Theorem $(\ref{mple})$ proves the existence of a solution to system $(\ref{exe})$ .
\end{pf}

\end{document}